\documentclass[11pt]{article}
\usepackage[letterpaper,hmargin=.95in,vmargin=.75in]{geometry}
\usepackage{amsmath,amssymb,amsthm,booktabs}
\usepackage[hidelinks]{hyperref}
\newtheorem{theorem}{Theorem}
\newtheorem{claim}{Claim}
\theoremstyle{definition}
\newtheorem{definition}{Definition}
\DeclareMathOperator{\area}{area}
\DeclareMathOperator{\dinv}{dinv}
\DeclareMathOperator{\sign}{sign}

\newcommand{\D}{\mathcal D}

\newcommand{\z}[1]{\#_0(#1)}
\newcommand{\ty}[1]{\mathsf{#1}}
\hypersetup{pdftitle={A combinatorial proof of area and dinv symmetry for Dyck paths}}
\begin{document}

\title{A combinatorial proof of area and dinv symmetry on Dyck paths}
\author{Graham Hawkes}
\date{}

\maketitle

\begin{abstract}
We construct an involution on Dyck paths that interchanges the statistics of area and dinv.
\end{abstract}

\section{Introduction}

\begin{definition}
A Dyck path of semilength \(n\) may be encoded via its area word~\cite{Haglund}. Precisely, define a \emph{Dyck word} \(W=w_1\cdots w_n\) to be a sequence of nonnegative integers with \(w_1=0\) and \(w_{i+1}\le w_i+1\) for \(1\le i<n\). Write \(\D_n\) for the set of Dyck words of length \(n\). Given any sequence of nonnegative integers \(W=w_1\cdots w_n\), define
\[
 \area(W)=\sum_{i=1}^{n}w_i,\qquad
 \dinv(W)=\#\{(i,j):1\le i<j\le n,\ w_i-w_j\in\{0,1\}\}.
\]
\end{definition}

\begin{theorem}\label{th:main}
For each \(n\ge0\) there is an involution \(\Phi_n\) on \(\D_n\) that interchanges \(\area\) and \(\dinv\).
\end{theorem}

\section{Proof}
\begin{definition}
Define an \emph{affine word} to be a sequence of nonnegative integers \(w_1\cdots w_n\) satisfying \(w_{i+1}\le w_i+1\) for \(1\le i<n\). Define a \emph{dual word} to be a sequence of nonnegative integers \(w_1\cdots w_n\) satisfying \(w_{i+1}>w_i+1\) for \(1\le i<n\)~\cite{Hawkes}. A \emph{positive word} in either case refers to one with no $0$s. 
\end{definition}

\begin{definition}\label{def:grammar}
Define a \emph{sentence} to be a list of punctuated words.  Each word may be Dyck, affine, or dual, may occur inside or outside of parentheses, and may be underlined or not.  We use the symbols \(\ty Y,\ty F,\ty U, \ty W,  (\ty W),\underline{\ty W}\) to refer to Dyck, affine, dual or arbitrary word types and their possible punctuations. The punctuation $\ty W^+$ indicates that only positive words of that type are allowed. 
\end{definition}

\begin{definition}
A word $W$ is equipped with a \emph{footnote}, $[j]$, where $j$ is a primed or unprimed nonnegative integer (write $W_{[j]}$). If not specified, the footnote is $[0]$. The \emph{value} of a footnote is its numerical value and its \emph{size}, indicated $\widehat{\cdot}$, is that value if it is unprimed and one more than it if primed.
\end{definition}

\begin{definition}
Define the \emph{type} of a sentence to be the ordered list of the types of its words, including their punctuation and footnotes. Let the \emph{content} of a word or sentence be the multiset of its entries. 
\end{definition}

\begin{definition}
Suppose $S$ is a sentence of fixed type and let $W$ be the
concatenation of its words. 
Across all footnotes in $S$ let $\sum q$ be the total sum of the values of footnotes on dual words, $\#q'$ the number of primed footnotes on dual words,  $\sum t$ the total sum of the footnote values on Dyck or affine words and $\#t'$ the number of primed footnotes on Dyck or affine words.  Define
\begin{align*}
\area(S) &= \area(W)
  + \#\text{parenthesized letters}
  + \sum {q} +\#q',\\
\dinv(S) &= \dinv(W)
  + \#\text{underlined zeros}
  + \sum t + \#q' ,\\
\sign(S) &= (-1)^d,\quad\text{where} \qquad
d = \#\text{letters in dual words}
  + \#q' + \#t'.
\end{align*}
\end{definition}

\begin{definition}
Define the length of a word  \(W\)  to be its number of entries, denoted $|W|$,  and its \emph{positive length} to be \(|W|-\z W\). The letter length of a sentence is the sum of the lengths of its words. 
\end{definition}

\begin{claim}\label{cl:cut}
Let $\ty T$ be a list of word types none of which are underlined. There is a bijection $\delta$ from
sentences of type
$
\ty Y,\ty T
$
with a nonempty first word to sentences of type
$
\underline{\ty Y},\underline{\ty  T},(\ty Y),
$
preserving the intervening words, their parentheses
and their footnotes, as well as sign, area and dinv, and decreasing
letter length by one.
\end{claim}

\begin{proof}
Write the first word as $0AD$, where $A$ is 
the maximal length positive affine word. 
Remove $0A$ from the first word and underline all resulting
words. Place the word obtained by lowering each
entry of $A$ by one in parentheses at the end of the sentence.
\end{proof}

\begin{definition}
Let \(X,Y\) be families whose elements are equipped with area, dinv and sign, and write \(X^+,X^-\) for their positive and negative parts, and similarly for \(Y\). Define a \emph{signed matching} \(X\simeq Y\) to be an area- and dinv-preserving bijection
\[
                     X^+\sqcup Y^-\longrightarrow X^-\sqcup Y^+.
\]
Define a \emph{cancellation} in \(X\) to be such a bijection from a subset of \(X^+\) to a subset of \(X^-\). 
\end{definition}

\begin{claim}\label{cl:cancel}
For each nonempty word \(Q\), the sentences of type \(\ty F,\ty U\) (resp. \(\ty U,\ty F\) ) with concatenation \(Q\) admit a sign-reversing involution $\beta$  (resp. $\beta'$) preserving the concatenation of the sentences.
\end{claim}
\begin{proof}
Move either the last letter of the first word onto the beginning of the second word or else the first letter of the second word onto the end of the first word. 
\end{proof}

\begin{definition}
Define a \emph{signed enlargement} from \(X\) to \(Y\) to be an embedding preserving sign, area and dinv together with a cancellation pairing every element outside its image. Pairing each element of \(X\) with its image and using this cancellation gives \(X\simeq Y\).
\end{definition}

\begin{claim}\label{cl:enlarge}
Let $\ty T$ be a list of word types.
There is a signed enlargement $\zeta$ from sentences of type
$\underline{\ty Y},\ty T$
to sentences of type
$\ty U^+,\underline{\ty F}, \ty T$
that fixes the words of types $\ty T$ and preserves
letter length. The same result holds if $\ty Y$ and $\ty F$ are not underlined above. 
\end{claim}

\begin{proof}
Define the embedding by prepending an initial empty word $\varnothing$ to the sentence. For sentences outside the image,
the concatenation of the first two words is nonempty and
begins with a positive letter. Thus $\beta'$ gives the
required cancellation.
\end{proof}

\begin{claim}\label{cl:alpha}
There is a bijection \(\alpha\) on sentences of type \(\ty U,\ty U\) that preserves content, sign, area and dinv and interchanges the full lengths of the two words. Cf. ~\cite[Cor.~3.19]{Hawkes}. 
\end{claim}
\begin{proof}
For each maximal consecutive string of integers of odd length greater than one, and of length one if that entry appears once in the content, switch which word has the odd entries from that string and which has the even ones. 
\end{proof}

\begin{claim}\label{cl:star}
There is a bijection \(\alpha_0\) on sentences of type \(\underline{\ty U},\ty U\) that preserves content, sign, area and dinv and interchanges the positive lengths of the two words.
\end{claim}
\begin{proof}
If both words or neither word contains a $0$, apply $\alpha$.  Otherwise ignore the $0$ and apply $\alpha$; if the $0$ now lies next to a $1$, move it to the beginning of the other word. 
\end{proof}

\begin{claim}\label{cl:mixed}
There is a bijection $\gamma$ from sentences of type
$\ty U,\ty F$ to sentences of type $\ty F,\ty U$ that
preserves content, sign, area and dinv and interchanges
the lengths of the two words.
\end{claim}

\begin{proof}
Fix nonnegative integers $r,s$ and a content of size $r+s$ and consider the two families
of sentences of type $\ty F,\ty U,\ty U,\ty F$ with respective
word lengths
$
(a,b,r,h)$ and 
$(a,r,b,h),$
where $a,b,h$ are nonnegative, $a+b+h=s$, and $b$ is even.
Applying $\alpha$ to the middle two words gives a bijection
between these families.
The composition $\beta'\circ\alpha\circ\beta$, with $\beta$
acting on the first two words and $\beta'$ on the last two (on which it must be able to act by the choice of $b$ even),
gives a bijection between the same families after omitting
sentences with word lengths
$
(0,0,r,s)$
and $
(s,r,0,0)$ from the first and second families respectively. 
Both bijections preserve content, area and dinv.
For each fixed content the families are finite, so
we may use Garsia--Milne to cancel the second bijection from the first
to give a bijection between the omitted families, whence the desired bijection can be constructed.
\end{proof}

\begin{claim}[Dual passage]\label{cl:dual}
Fix nonnegative integers $m,s$.
There is a content-preserving signed matching $\phi^*$ from sentences of type
$
                       (\underline{\ty U}),\ty Y
$
with word lengths $(s,m)$ to the disjoint union of sentences of types
\[
\begin{array}{ll}
(\ty Y,(\ty U)_{[j]}) & 0\le j\le s,\\[3pt]
(\ty Y,(\ty U)_{[(j-1)']}) & 1\le j\le s,
\end{array}
\]
with word lengths $(m+j,s-j)$ in both cases.
\end{claim}
\begin{proof}
Fix a content of size $m+s$.
Define $L$ and $R$ to be the families of sentences with this content and types
$(\underline{\ty U}_{[s]},\ty U,\ty F)$
and
$(\underline{\ty U},\ty F,\ty U_{[s]})$,
respectively,
such that, in either case, $h$ is the length of the affine word and $b$ is the positive length of the dual word immediately preceding it  and $b+h\ge m$.
The composition $\gamma \circ \alpha_0$ (applied to the leftmost two words and then the rightmost two words), along with moving the footnote from the first word to the last, gives a content, sign, area, and dinv preserving bijection $L\to R$. 

 Define the bijection
$
f:L^+\sqcup R^-\longrightarrow L^-\sqcup R^+
$
to be this composition on $L^+$ and its inverse on $R^-$.
Let $X\subseteq L$ consist of the sentences whose second word
is empty and whose third word is a Dyck word of length $m$ (including the empty Dyck word when $m$=0).
Let $Y\subseteq R$ consist of the sentences whose first two
words concatenate to a Dyck word.
Define
\[
g:(L^+\setminus X^+)\sqcup(R^-\setminus Y^-)
\longrightarrow
(L^-\setminus X^-)\sqcup(R^+\setminus Y^+)
\]
by applying $\beta'$ to the last two words of sentences in
$L^+\setminus X^+$ and to the first two words of sentences
in $R^-\setminus Y^-$.
Note the only way for $g$ to hypothetically leave $L \setminus X$ is if it moves a $0$ right (could create a Dyck sequence, but since $b+h \ge m$ it would be longer than $m$) or to the left (but the  affine word must still have length at least $m$ since $g$ does not act on $X$). The only way for $g$ to leave $R \setminus Y$ would be for $0$ to move between factors, but this is impossible. 
Use Garsia--Milne to cancel $g$ against $f$, obtaining
a bijection
$
X^+\sqcup Y^-\longrightarrow X^-\sqcup Y^+,
$
and hence a signed matching between $X$ and $Y$.

Deleting the empty second word in $X$, parenthesizing its
first word and removing its footnote gives
the family of sentences with type $
                       (\underline{\ty U}),\ty Y
$
with word lengths $(s,m)$.   Concatenating the first two words of $Y$
into a Dyck word
of length $m+j$, parenthesizing the remaining dual word of
length $s-j$, and replacing its footnote by $[j]$ or
$[(j-1)']$,  if the original first word was
empty or $0$, respectively, gives the second family in the claim. 
Both operations preserve sign, area, dinv, and content. 
\end{proof}

\begin{claim}[Affine passage]\label{cl:aff}
Fix nonnegative integers $m,s$.
There is a content-preserving signed matching $\phi$ from
sentences of type $\underline{\ty F},\ty Y$ with word lengths
$(s,m)$ to the disjoint union of sentences of types
\[
\begin{array}{ll}
(\ty Y,\ty F_{[j]}) & 0\le j\le s,\\[3pt]
(\ty Y,\ty F_{[(j-1)']}) & 1\le j\le s,
\end{array}
\]
with word lengths $(m+j,s-j)$ in both cases.
\end{claim}

\begin{proof}
For $0\le d\le s$, fix a content of size $m+d$. 
Define $L$ and $R$ to be the families of sentences with this content and types
$(\underline{\ty U},\ty U,\ty F_{[p]})$
and
$(\underline{\ty U},\ty F_{[p]},\ty U)$,
such that in either case $h$ is the affine word's length, $b$ is the positive length of the dual word immediately
preceding it and  $p=b+h-m\ge0$.  The composition $\gamma\circ\alpha_0$ (applied
first to the leftmost two words and then to the rightmost two words and keeping the footnote on the affine word)
gives a content, sign, area, and dinv preserving bijection $L\to R$.

Define
$f:L^+\sqcup R^-\to L^-\sqcup R^+$
to be this composition on $L^+$ and its inverse on $R^-$.
Let $X\subseteq L$ consist of sentences whose last two words
concatenate to a Dyck word, and let $Y\subseteq R$ consist
of sentences whose first word is empty and whose second
word is Dyck of length $m$.
Define
\[
g:(L^+\setminus X^+)\sqcup(R^-\setminus Y^-)
\longrightarrow
(L^-\setminus X^-)\sqcup(R^+\setminus Y^+)
\]
by applying $\beta'$ to the last two words in $L$ and
the first two in $R$ (and in the latter case increasing (resp. decreasing) the footnote entry on the affine word by $1$ if a $0$ is passed right (resp. left)). This is well-defined by the same logic as in the previous proof.
Now, use Garsia--Milne to cancel $g$ against $f$, obtaining
a bijection
$X^+\sqcup Y^-\longrightarrow X^-\sqcup Y^+$,
and hence a signed matching between $X$ and $Y$.

Deleting the empty first word in $Y$ 
 gives the family $\mathcal{A}$ of sentences of type
$\ty Y,\ty U$ with word lengths $(m,d)$.
Concatenating the last two words in $X$ into a Dyck word
of length $m+j$ and giving the resulting Dyck word the footnote $[j]$
if the original second word was empty or $[(j-1)']$ if it was $0$
gives a family  $\mathcal{B}$ of sentences of type
\[
\begin{array}{ll}
(\underline{\ty U},\ty Y_{[j]}) & 0\le j\le d,\\[3pt]
(\underline{\ty U},\ty Y_{[(j-1)']}) & 1\le j\le d,
\end{array}
\]
with word lengths $(d-j,m+j)$ in both cases.  The signed matching from $X$ to $Y$ thus gives a signed matching $\pi$ from $\mathcal A$ to $\mathcal B$.

Let $P$ and $Q$ be the families of sentences of types
$\underline{\ty F}:\mathcal{A}:\ty F$ and
$\underline{\ty F}:\mathcal{B}:\ty F$, respectively,
with outer words of any lengths $a,c$, such that $a+c=s-d$, for some $0 \le d \le s$
and with a fixed content of size $m+s$.
Applying $\pi$ to the middle two words gives a bijection
\[
f_1:P^+\sqcup Q^-\longrightarrow P^-\sqcup Q^+.
\]

Let $X_1\subseteq P$ consist of the sentences whose last
two words are empty, and let $Y_1\subseteq Q$ consist
of the sentences whose first two words are empty.
Define
\[
g_1:(P^+\setminus X_1^+)\sqcup(Q^-\setminus Y_1^-)
\longrightarrow
(P^-\setminus X_1^-)\sqcup(Q^+\setminus Y_1^+)
\]
by applying $\beta'$ to the last two words of sentences in
$P$  and $\beta$ to the first two words of sentences in
$Q$. Note these two pairs of words are both not underlined or both underlined respectively so dinv is preserved.
We may now use Garsia--Milne to cancel $g_1$ against $f_1$,
obtaining a signed matching between $X_1$ and $Y_1$, which implies one between the families in the claim, after transferring the footnote to the second word for the second class.
\end{proof}

\begin{definition}
For a list $\boldsymbol a=(a_1,\ldots,a_k)$ of nonnegative integers,
define $\D_{\boldsymbol a}$ to consist of sentences
$(D,A_2,\ldots,A_k,B_k,\ldots,B_2)$ of types
$\ty Y,\ty F,\ldots,\ty F,(\ty U),\ldots,(\ty U)$,
with $|D|=a_1$ and $|A_i|+|B_i|=a_i$ for $i>1$.
\end{definition}

\begin{definition}\label{def:archive}
A \emph{record} $r$ is a list of primed or unprimed entries of even length, $[t_1,\ldots,t_k, q_k,\ldots, q_1]$ with last entry unprimed.
It acts \emph{validly} on a list $\boldsymbol a$ of nonnegative integers when, after removing leading $0$s, the resulting list $(a_1,\ldots,a_k)$ has half the length of the record and the following \emph{output} list has no negative values,
\begin{equation}\label{eq:row}
 \theta(r,\boldsymbol a)=\left(\sum (\widehat q_i+\widehat t_i),\ a_1-1-(\widehat q_1+\widehat t_1),\ a_2-(\widehat q_2+\widehat t_2),\ldots,a_k-(\widehat q_k+\widehat t_k)\right).
\end{equation}
An \emph{archive} of length $n$ is a list of $n$ records  such that  each successive record  acts validly on the output list produced by the last record, starting with the list $(n)$.
Let $\mathcal R_n$ be the family of valid archives of length $n$, where the following statistics are taken across all records in the archive
\[
 \area(R)=\sum q+\#q',\qquad \dinv(R)=\sum t+\#q',\qquad
 \sign(R)=(-1)^{\#q'+\#t'}.
\]
\end{definition}

\begin{claim}\label{cl:short}
For any list $\boldsymbol a$ with positive sum, there is a signed matching $\sigma$ from $\D_{\boldsymbol a}$ to pairs $(r,S)$, where $r$ is a valid record on $\boldsymbol a$ and $S\in\D_{\theta(r,\boldsymbol a)}$ (statistics on pairs are computed the natural way).
\end{claim}

\begin{proof}
First, we may assume $a_1 \neq 0$.  Indeed, if $a_1=0$ we can easily construct a signed matching with $\D_{(a_2,\ldots,a_k)}$: Starting with $\D_{\boldsymbol a}$ replace it with the family obtained by acting elementwise on its sentences by deleting the first empty word, removing the last dual word, unparenthesizing it, adding one to each of its letters and moving it to the front of the sentence; and then applying Claim \ref{cl:enlarge} to the resulting family. Repeating inductively shows our assumption is valid. 

Now apply $\delta$ to $\D_{\boldsymbol a}$, then enlarge the resulting family by acting on the first word with $\zeta$.
Next, apply the operation of removing the first word (type $ \ty U^+$) lowering its entries by one and placing it in parentheses at the end of the sentence. Now, initialize a record $r$ of length $2k$, remove the parentheses around second to last word (sending type $(\ty Y) \rightarrow \ty Y$) of length say $q_1$ and enter $q_1$ as the last entry of $r$.
Then pass this Dyck word left over $k-1$ dual words using $\phi^*$ and then over $k$ affine words using $\phi$ so that in each successive family the Dyck word lies one position farther left. At each step pull off the footnote that was left on the word that was passed and add it to $r$, populating the record right to left.

It remains to verify \eqref{eq:row}.  Before the applications of  $\phi^*$ and $\phi$ the moving Dyck word had size $q_1$, the sum of the lengths of the first and last word (which were created from the original first word after splitting off a $0$ and a word of length $q_1$) was $(a_1-1-q_1)$ and all other words had their original sizes.  After the application of these maps, the length of the moving Dyck word grew by the sizes of the first $2k-1$ record entries, the (now second) word shrank by the size of the first record entry, and the  other sums of relevant pairs of lengths shrank by the sums of the sizes of the corresponding record entries. This verifies \eqref{eq:row}.  Since these lengths are nonnegative, $r$ is valid and $S\in\D_{\theta(r,\boldsymbol a)}$. 

Conversely, to see that the inverse operation is well-defined note that since $S\in\mathcal D_{\theta(r,\boldsymbol a)}$,
the first coordinate of $\theta(r,\boldsymbol a)$ guarantees
that the moving Dyck word has sufficient length to lose as it passes rightward, and all of the coordinates together imply that the resulting sentences have the word lengths needed to belong to $\mathcal D_{\boldsymbol a}$.
\end{proof}

\begin{proof}[Proof of Theorem~\ref{th:main}]
Start with $\D_{(n)}$ and an empty archive. Repeatedly apply $\sigma$, archiving the record each time. Each application lowers letter length by one, so after $n$ steps only empty words remain. Their composition gives a signed matching from $\D_n$ to $\mathcal R_n$, that is, a bijection
$
 f:\D_n\sqcup\mathcal R_n^-\longrightarrow\mathcal R_n^+
$
preserving area and dinv.  For $R\in\mathcal R_n$, define $\tau(R)$ by reversing the order of the numerical values in each record, keeping the primes in place. This fixes each sum $\widehat q_i+\widehat t_i$ so equation \eqref{eq:row} implies $\tau$ preserves validity. By construction $\tau$ is a sign-preserving involution on $\mathcal R_n$ that interchanges area and dinv, and $f^{-1}\tau f$ is an involution on $\D_n\sqcup\mathcal R_n^-$ that also interchanges area and dinv.  Thus using Garsia--Milne to cancel $\tau$ on $\mathcal R_n^-$ against $f^{-1}\tau f$ gives an area and dinv interchanging bijection $\Phi_n$ on $\D_n$. Since $\tau$ and $f^{-1}\tau f$ are both involutions, the inverse cancellation uses the same two maps. Hence $\Phi_n^{-1}=\Phi_n$.

\end{proof}

\section*{Acknowledgements}
The author used Codex to write Python scripts to check combinatorial equalities in finite cases, to double-check certain mathematical arguments before committing to developing them fully, and to help with spelling, grammar, and \LaTeX{} formatting.  In addition, the original version of this proof relied on results from \cite{Hawkes} and \cite{Hawkes1},  both of which used AI tools to varying extents (see acknowledgments there); these dependencies have been fully removed in the present version and replaced with much weaker statements that are proved directly here.  Precisely, Claims~\ref{cl:alpha} and~\ref{cl:cancel} replace earlier dependencies on  \cite{Hawkes} and \cite{Hawkes1}, respectively.

\end{document}